\documentclass{amsart}

\usepackage{amsmath,amssymb,amscd,amsfonts}

\newtheorem{theorem}{Theorem}
\newcommand{\bt}{\begin{theorem}}
\newcommand{\et}{\end{theorem}}
\newtheorem{lemma}{Lemma}
\newcommand{\bl}{\begin{lemma}}
\newcommand{\el}{\end{lemma}}
\newtheorem{corollary}{Corollary}
\newcommand{\bc}{\begin{corollary}}
\newcommand{\ec}{\end{corollary}}
\newtheorem{problem}{Problem}
\newcommand{\bprob}{\begin{problem}}
\newcommand{\eprob}{\end{problem}}

\newtheorem{example}{Example}
\newcommand{\bex}{\begin{example}}
\newcommand{\eex}{\end{example}}

\newcommand{\beq}{\begin{equation}}
\newcommand{\eeq}{\end{equation}}
\newcommand{\benum}{\begin{enumerate}}
\newcommand{\eenum}{\end{enumerate}}
\newcommand{\N}{\ensuremath{ \mathbf N }}

\newcommand{\R}{\ensuremath{\mathbf R}}

\newcommand{\mce}{\ensuremath{ \mathcal E}}

\newcommand{\mba}{\ensuremath{ \mathbf a}}
\newcommand{\mbb}{\ensuremath{ \mathbf b}}
\newcommand{\mbc}{\ensuremath{ \mathbf c}}

\newcommand{\mbe}{\ensuremath{ \mathbf e}}

\newcommand{\mbf}{\ensuremath{ \mathbf f}}
\newcommand{\mbv}{\ensuremath{ \mathbf v}}
\newcommand{\mbw}{\ensuremath{ \mathbf w}}
\newcommand{\mbx}{\ensuremath{ \mathbf x}}
\newcommand{\mby}{\ensuremath{ \mathbf y}}

\newcommand{\Rn}{\ensuremath{ \mathbf{R}^n }}
\DeclareMathOperator{\colsum}{\text{colsum}}
\DeclareMathOperator{\rowsum}{\text{rowsum}}
\newcommand{\bsmallmat}{\left(\begin{smallmatrix}}
\newcommand{\esmallmat}{\end{smallmatrix}\right)}

\DeclareMathOperator{\card}{\text{card}}

\newcommand{\bmat}{\left(\begin{matrix}}
\newcommand{\emat}{\end{matrix}\right)}

\DeclareMathOperator{\qqand}{\qquad\text{and}\qquad}

\DeclareMathOperator{\vectoran}{\left( \begin{matrix} a_1 \\ \vdots \\ a_n \end{matrix}\right)}
\DeclareMathOperator{\vectorbn}{\left( \begin{matrix} b_1 \\ \vdots \\ b_n \end{matrix}\right)}
\DeclareMathOperator{\vectorcn}{\left( \begin{matrix} c_1 \\ \vdots \\ c_n \end{matrix}\right)}

\DeclareMathOperator{\vectorsmallan}{\left( \begin{smallmatrix} a_1 \\ \vdots \\ a_n \end{smallmatrix}\right)}
\DeclareMathOperator{\vectorsmallbn}{\left( \begin{smallmatrix} b_1 \\ \vdots \\ b_n \end{smallmatrix}\right)}
\DeclareMathOperator{\vectorsmallcn}{\left( \begin{smallmatrix} c_1 \\ \vdots \\ c_n \end{smallmatrix}\right)}

\DeclareMathOperator{\vectorsmallyn}{\left( \begin{smallmatrix} y_1 \\ \vdots \\ y_n \end{smallmatrix}\right)}

\newcommand{\mbo}{\ensuremath{\mathbf 0}}
\DeclareMathOperator{\col}{\text{col}}

\title{The Muirhead-Rado inequality, 1: Vector majorization and the permutohedron}
\author{Melvyn B. Nathanson}
\address{Lehman College (CUNY), Bronx, New York 10468} 
\email{melvyn.nathanson@lehman.cuny.edu}



\thanks{Supported in part by PSC CUNY Grant \# 66197-00 54.}

\begin{document}

\begin{abstract}
Let $\mathbf{a}$  and $\mathbf{b}$ be vectors in $\mathbf{R}^n$ with nonnegative coordinates.  
Permuting the coordinates, one can assume that $a_1 \geq \cdots \geq a_n$ 
and $b_1 \geq \cdots \geq b_n$.  The vector $\mathbf{a}$ majorizes the vector $\mathbf{b}$  
 if $\sum_{i=1}^k b_i \leq \sum_{i=1}^k a_i$ for all $k \in \{1,\ldots,n-1\}$ 
 and  $\sum_{i=1}^n b_i = \sum_{i=1}^n a_i$. 
This paper exposes theorems of 
Hardy-Littlewood-P\'olya and of Rado: The following are equivalent: 
(1) The vector $\mathbf{a}$ majorizes the vector $\mathbf{b}$, 
(2)  $P\mathbf{a} = \mathbf{b}$ for some 
doubly stochastic matrix $P$, 
and (3) $ \mathbf{b}$ is in the 
$S_n$-permutohedron generated by $\mathbf{a}$
\end{abstract}

\maketitle

\section{Permutation matrices and the permutohedron}

This paper describes majorization results of 
Hardy-Littlewood-P\'olya~\cite{hard-litt-poly29,hard-litt-poly88}  
and Rado~\cite{rado52}  
that will be applied in~\cite{nath24} to prove beautiful inequalitites of Muirhead~\cite{muir03}  
and Rado~\cite{rado52}  
that generalize the arithmetic and 
geometric mean  inequality.    
An encyclopedic reference is Marshall-Olkin-Arnold~\cite {mars-olki-arno-11}. 

A matrix  $P = \bmat p_{i,j} \emat$ with real coordinates is 
\emph{nonnegative}\index{matrix!nonnegative} 
if $p_{i,j} \geq 0$ for all $i$ and $j$. 
A nonnegative $n\times n$ matrix $P = \bmat p_{i,j} \emat$ is 
\emph{column stochastic}\index{column stochastic} if its column sums satisfy 
\[
\colsum_j(P) = \sum_{i=1}^n p_{i,j} = 1 \qquad \text{for all $j \in \{1,\ldots, n\}$.}
\]
A nonnegative $n\times n$ matrix $P = \bmat p_{i,j} \emat$ is 
\emph{row stochastic}\index{row stochastic} if its row sums satisfy 
\[
\rowsum_i(P) = \sum_{j=1}^n p_{i,j} = 1 \qquad \text{for all $i \in \{1,\ldots, n\}$}.
\]
A nonnegative $n\times n$ matrix $P = \bmat p_{i,j} \emat$ is 
\emph{doubly stochastic}\index{doubly stochastic} if  
it is both column stochastic and row stochastic.  

We begin with some simple lemmas.  

\bl                \label{Muirhead:lemma:DoublyStochasticProduct}
Let $P$ and $Q$ be $n \times n$ matrices.  
If $P$ and $Q$ are column stochastic, then $PQ$ is column stochastic. 
If $P$ and $Q$ are row stochastic, then $PQ$ is row stochastic. 
If $P$ and $Q$ are doubly stochastic, then $PQ$ is doubly stochastic. 
\el

\begin{proof}
Let $P = \bmat p_{i,j} \emat$ and $Q = \bmat q_{i,j} \emat$. 
If $P$ and $Q$ are column stochastic  $n \times n$ matrices, then for all $j \in \{1,\ldots, n\}$ we have 
\begin{align*}
\colsum_j(PQ)&  = \sum_{i=1}^n (PQ)_{i,j} = \sum_{i=1}^n \sum_{k=1}^n p_{i,k}q_{k,j} \\
& =  \sum_{k =1}^n  q_{k,j} \sum_{i=1}^n p_{i,k}  =  \sum_{k =1}^n q_{k,j} \colsum_k(P) \\
& =  \sum_{k =1}^n q_{k,j} = \colsum_j(Q) = 1
\end{align*}
and so $PQ$ is column stochastic.

Similarly, if $P$ and $Q$ are row stochastic  $n \times n$ matrices, 
then for all $i \in \{1,\ldots, n\}$ we have
\begin{align*}
\rowsum_i(PQ)&  = \sum_{j=1}^n (PQ)_{i,j} = \sum_{j=1}^n \sum_{k=1}^n p_{i,k}q_{k,j} \\
& =  \sum_{k =1}^n p_{i,k}  \sum_{j=1}^n q_{k,j} =  \sum_{k =1}^n  p_{i,k} \rowsum_k(Q)  \\
& =  \sum_{k =1}^n  p_{i,k} = \rowsum_i(P) = 1
\end{align*}
and so $PQ$ is row stochastic. 

It follows that $PQ$ is doubly stochastic if $P$ and $Q$ are doubly stochastic.  
This completes the proof. 
\end{proof}

Let  $\mce = \{\mbe_1,\ldots, \mbe_n\}$ be the standard basis for $V = \Rn$.    
For every permutation  $\sigma$  in the symmetric group $S_n$, define the linear transformation 
$P_{\sigma}:V \rightarrow V$ by 
\[
P_{\sigma}(\mbe_j) = \mbe_{\sigma(j)}
\]
for all $j \in \{1,\ldots, n\}$. 
If $\mbx = \sum_{j=1}^n x_j\mbe_j \in V$, then 
\begin{align*}
P_{\sigma}(\mbx) & = P_{\sigma}\left( \sum_{j=1}^n x_j\mbe_j \right) 
=   \sum_{j=1}^n x_j P_{\sigma}\left( \mbe_j \right) \\ 
& =   \sum_{j=1}^n x_j \mbe_{\sigma(j)} 
=   \sum_{j=1}^n x_{\sigma^{-1}(j)} \mbe_j.  
\end{align*}
Equivalently,  
\[
P_{\sigma}  \left( \begin{matrix} x_1 \\ \vdots \\ x_n \end{matrix}\right)    = \bmat x_{\sigma^{-1}(1)} \\ \vdots \\ x_{\sigma^{-1}(n)} \emat.
\] 

We also denote by $P_{\sigma} $ the matrix of this linear transformation 
with respect to the standard basis \mce. 
The  $n \times n$  matrix $P_{\sigma} $ is called a 
\emph{permutation matrix}\index{matrix!permutation}.    
Its $(i,j)$th coordinate is      
\[
\left(P_{\sigma}\right)_{i,j}  = \delta_{i,\sigma(j)} = \begin{cases}
1 & \text{if $i = \sigma(j)$}\\
0 & \text{if $i \neq \sigma(j)$}
\end{cases}
\]
where $\delta_{i,j}$ is the Kronecker delta.  
For all $j \in \{1,\ldots, n\}$ we have 
\[
\colsum_j(P_{\sigma}) = \sum_{i=1}^n \left(P_{\sigma}\right)_{i,j}  = \sum_{i =1}^n  \delta_{i,\sigma(j)} = 1. 
\]
For all $i \in \{1,\ldots, n\}$ we have 
\[
\rowsum_i(P_{\sigma}) = \sum_{j=1}^n \left(P_{\sigma}\right)_{i,j}  = \sum_{j=1}^n  \delta_{i,\sigma(j)} = 1. 
\] 
Thus, the permutation matrix $P_{\sigma}$ is doubly stochastic for all $\sigma \in S_n$. 

For example, if $\sigma = (1,2,3) \in S_3$, then 
\[
P_{\sigma} = \bmat 
0 & 0 & 1 \\ 
1 & 0 & 0 \\ 
0 & 1 & 0 \emat.
\]

\bl
For all permutations $\sigma, \tau \in S_n$, 
\[
P_{\sigma} P_{\tau} = P_{\sigma\tau}.  
\]
\el

\begin{proof}
For all $j \in \{1,\ldots, n\}$ we have 
\begin{align*}
\left( P_{\sigma} P_{\tau} \right)(\mbe_j) 
& =  P_{\sigma} \left( P_{\tau} (\mbe_j) \right) 
= P_{\sigma} \left( \mbe_{\tau(j)} \right) = \mbe_{\sigma( \tau(j))}  \\ 
& =  \mbe_{ (\sigma \tau) (j)}  = P_{\sigma\tau}(\mbe_j).
\end{align*}
Equivalently, for all $i,j \in \{1,\ldots, n\}$ we have 
\[
\left(P_{\sigma} P_{\tau}\right)_{i,j} = \sum_{k=1}^n \left(P_{\sigma} \right)_{i,k}\left(P_{\tau}\right)_{k,j} 
= \sum_{k=1}^n \delta_{i,\sigma(k)} \delta_{k,\tau(j)}  
= \delta_{i,\sigma \tau(j)}  = \left(P_{\sigma \tau}\right)_{i,j}. 
\]
This completes the proof.  
\end{proof}

\bl       \label{Muirhead:lemma:sigma-tau}
Let $A$ be a $m \times n$ matrix.  
Let $\sigma$ be a permutation in $S_n$ and let $\tau$ be a permutation in $S_m$.  
\benum
\item[(a)]
The $j$th column of the $m \times n$ matrix $AP_{\sigma}$ is the 
$\sigma(j)$th column of the matrix $A$.  
\item[(b)]
 The $i$th row of the $m \times n$ matrix $P_{\tau}A$ is the 
$\tau^{-1}(i)$th row of the matrix $A$. 
\item[(c)]
The $(i,j)$th coordinate of the matrix $P_{\tau}AP_{\sigma}$ is $a_{\tau^{-1}(i),\sigma(j)}$. 
\eenum
\el

\begin{proof}
Matrix multiplication.  
\end{proof}

The classical \emph{permutohedron}\index{permutohedron}\index{permutohedron!classical permutohedron}
 is the convex hull of the set 
of $n!$ lattice points in \Rn\ whose coordinates are the permutations of the set $\{1,2,3,\ldots, n\}$, 
that is, the convex hull of the set 
\[
\left\{ P_{\sigma}\bmat 1 \\ 2\\ \vdots\\ n \emat: \sigma \in S_n \right\} 
= \left\{ \bmat {\sigma^{-1}(1)} \\ {\sigma^{-1}(2) } \\ \vdots\\ {\sigma^{-1}(n) } \emat : \sigma \in S_n \right\} \\
= \left\{ \bmat {\sigma(1)} \\ {\sigma(2) } \\ \vdots\\ {\sigma(n) } \emat : \sigma \in S_n \right\}. 
\]
Every vector in this set lies on the hyperplane 
\[
x_1 + x_2 + \cdots + x_n = 1 + 2 + \cdots + n = \frac{n(n+1)}{2}.
\]

For example,  the classical permutohedron in $\R^2$ is the convex hull of the set 
\[
\left\{ P_{\sigma}\bmat 1 \\ 2 \emat: \sigma \in S_2 \right\} 
\]
that is, the line segment 
\[
\left\{ (1-t)  \bmat 1 \\ 2  \emat + t  \bmat 2 \\ 1 \emat : 0 \leq t \leq 1 \right\} 
 = \left\{  \bmat  1+t \\ 2 - t \emat : 0 \leq t \leq 1 \right\} . 
\]
The classical permutohedron in $\R^3$ is the convex hull of the set 
\[
\left\{ P_{\sigma}\bmat 1 \\ 2\\ 3 \emat: \sigma \in S_3 \right\} 
= \left\{  \bmat 1 \\ 2\\ 3 \emat, \bmat 1 \\ 3\\ 2 \emat, \bmat 2 \\ 1 \\ 3 \emat, \bmat  2\\ 3 \\ 1 \emat, 
\bmat 3 \\1 \\ 2 \emat, \bmat 3 \\ 2 \\ 1\emat \right\}.
\]

Let $\mba = \vectorsmallan \in \Rn$.  
The \emph{permutohedron} $K(\mba)$ 
constructed from the vector \mba\ 
is the convex hull of the set of vectors in \Rn\ 
whose coordinates are the permutations 
of the set $\{a_1, a_2, a_3,\ldots, a_n\}$, 
that is, the convex hull of the set 
\[
\left\{ P_{\sigma}\bmat a_1 \\ a_2\\ \vdots\\ a_n \emat: \sigma \in S_n \right\} 
= \left\{ \bmat a_{\sigma(1)} \\ a_{\sigma(2) } \\ \vdots\\ a_{\sigma(n) } \emat 
: \sigma \in S_n \right\}.
\]
Every vector in this set lies on the hyperplane 
\[
H = \left\{  \left( \begin{matrix} x_1 \\ \vdots \\ x_n \end{matrix}\right)    \in \Rn : x_1 + x_2 + \cdots + x_n = a_1 + a_2 + \cdots + a_n \right\}. 
\]

Let $G$ be a subgroup of $S_n$. 
The \emph{$G$-permutohedron}\index{permutohedron!$G$-permutohedron} 
constructed from the vector $\mba = \vectorsmallan$, denoted $K_G(\mba)$,  
is the convex hull of the set of vectors in \Rn\ whose coordinates are the $G$-permutations 
of the vector \mba, 
that is, the set 
\[
\left\{ P_{\sigma}\bmat a_1 \\ a_2\\ \vdots\\ a_n \emat: \sigma \in G \right\} 
=  \left\{ \bmat a_{\sigma(1)} \\ a_{\sigma(2) } \\ \vdots\\ a_{\sigma(n) } \emat 
: \sigma \in G \right\}.
\]
The $G$-permutohedron $K_G(\mba)$ is a subset of the permutohedron $K(\mba)$.

\section{Decreasing and increasing  vectors}  
For every permutation $\sigma \in S_n$ we have the permutation matrix $P_{\sigma}$ and 
for every vector $\mba = \vectorsmallan \in \R^n$ we have the vector 
with reordered    coordinates 
\[ 
P_{\sigma}(\mba) = \bmat a_{\sigma^{-1}(1)} \\ \vdots \\ a_{\sigma^{-1}(n)} \emat \in \R^n.  
\]
The  vector $P_{\sigma}(\mba)$ is called a 
\emph{rearrangement}\index{rearrangement} of the vector in \mba.  
Of special interest are vectors with decreasing coordinates 
and  vectors with increasing coordinates. 

Let $\mbx = \left( \begin{smallmatrix} x_1 \\ \vdots \\ x_n \end{smallmatrix}\right) \in \Rn$.  
The vector $\mbx$ is \emph{decreasing}\index{vector!decreasing} if 
\[
x_1 \geq x_2 \geq  \cdots \geq x_n.
\]
Associated to every vector $\mba = \vectorsmallan \in \Rn$ is a unique decreasing vector 
$\mba^{\downarrow} = \bsmallmat a_1^{\downarrow}  \\ \vdots \\ a_n^{\downarrow}\esmallmat \in \Rn $ 
obtained from \mba\ by a rearrangement of  coordinates. 
For every vector $\mba  \in \Rn$ there is a permutation $\sigma \in S_n$ 
such that $P_{\sigma}(\mba) = \mba^{\downarrow}$.
This permutation is unique if and only if the vector \mba\ has distinct coordinates.

The vector $\mbx = \left( \begin{smallmatrix} x_1 \\ \vdots \\ x_n \end{smallmatrix}\right)\in \Rn$ is \emph{increasing}\index{vector!increasing}
if 
\[
x_1 \leq x_2 \leq  \cdots \leq x_n.
\]
Associated to every vector 
$\mba = \vectorsmallan \in \Rn$ 
is a unique increasing vector 
$\mba^{\uparrow} = \bsmallmat a_1^{\uparrow}  \\ \vdots \\ a_n^{\uparrow}\esmallmat \in \Rn $ 
obtained from \mba\ by a rearrangement of  coordinates. 
For every vector $\mba   \in \Rn$ there is a permutation $\tau \in S_n$ 
such that $P_{\tau}\mba = \mba^{\uparrow}$.  
This permutation is unique if and only if the vector \mba\ has distinct coordinates.  

For example,   the  decreasing vector    obtained from the vector 
\[
\mba = \bmat a_1 \\ a_2 \\ a_3 \\ a_4 \emat = \bmat 3 \\ 1 \\ 2 \\ 4 \emat \in \R^4. 
\]
is 
\[
\mba^{\downarrow} = \bmat 4 \\ 3 \\ 2 \\ 1 \emat = \bmat a_4 \\ a_1 \\ a_3 \\ a_2 \emat 
= \bmat a_{\sigma^{-1}(1)} \\  a_{\sigma^{-1}(2)} \\ a_{\sigma^{-1}(3)} \\ a_{\sigma^{-1}(n)} \emat 
\]
where the rearrangement permutation $\sigma \in S_4$ is 
\[
\sigma^{-1} = \bmat 1  & 2 & 3 & 4 \\
4 & 1 & 3 & 2 \emat
\qqand 
\sigma = \bmat 1  & 2 & 3 & 4 \\
 2 & 4 &3 &  1 \emat.
\]
Thus,  $\mba^{\downarrow} = P_{\sigma}(\mba) $ and $\sigma$ is the 3-cycle $(1,2,4)$.

The increasing vector   obtained from \mba\ is 
\[
\mba^{\uparrow} = \bmat 1 \\ 2 \\ 3 \\ 4  \emat = \bmat a_2 \\ a_3 \\ a_1 \\ a_4 \emat 
= \bmat a_{\tau^{-1}(1)} \\  a_{\tau^{-1}(2)} \\ a_{\tau^{-1}(3)} \\ a_{\tau^{-1}(n)} \emat 
\]
where the rearrangement permutation $\tau \in S_4$ is 
\[
\tau^{-1} = \bmat 1  & 2 & 3 & 4 \\
2 & 3 & 1 & 4 \emat
\qqand 
\tau = \bmat 1  & 2 & 3 & 4 \\
3 &  1 &  2 & 4 \emat. 
\]
Thus, $\mba^{\uparrow} = P_{\tau}(\mba)$ and $\tau$ is the 3-cycle $(1,3, 2)$.

If $\mbb = \bsmallmat 2 \\1 \\ 3\\ 2 \esmallmat$, then 
\[
\mbb^{\downarrow}  = \bsmallmat  3 \\ 2 \\ 2 \\ 1 \esmallmat = P_{\sigma_1}(\mbb)   = P_{\sigma_2}(\mbb)
\]  
where $\sigma_1 = (1,2,4,3)$  and  $\sigma_2 = (1,3)(2,4)$, and    
\[
\mbb^{\uparrow}  = \bsmallmat  1 \\ 2 \\ 2 \\ 3 \esmallmat = P_{\tau_1}(\mbb)   = P_{\tau_2}(\mbb) 
\]
where $\tau_1 = (1,3,4, 2)$  and  $\tau_2 = (1,2)(3,4)$.  
The coordinates of \mbb\ are not distinct and the rearrangement permutations 
are not unique.  

Note that for all vectors $\mba, \mbb \in \R^n$ and all permutations $\sigma \in S_n$ 
we have 
\[
\left(\mba^{\downarrow} + \mbb^{\downarrow}\right)^{\downarrow} 
= \mba^{\downarrow} + \mbb^{\downarrow} 
 \qqand 
 \left(\mba^{\uparrow} + \mbb^{\uparrow}\right)^{\uparrow} 
= \mba^{\uparrow} + \mbb^{\uparrow}
\]
and 
\[
\left(P_{\sigma}(\mba)\right)^{\downarrow} = \mba^{\downarrow}  \qqand 
\left(P_{\sigma}(\mba)\right)^{\uparrow} = \mba^{\uparrow}.
\]

Let $\mbx = \left( \begin{smallmatrix} x_1 \\ \vdots \\ x_n \end{smallmatrix}\right) \in \Rn$  
and let $\tau  = (k,\ell) \in S_n$ be a  transposition with $k < \ell$.  
The  transposition $\tau$ \emph{increases the vector} $\mbx$ if $x_k > x_{\ell}$ and  
$P_{\tau}(\mbx)_k < P_{\tau}(\mbx)_{\ell}$. 
The  transposition $\tau$ \emph{decreases the vector} $\mbx$ if $x_k < x_{\ell}$ and  
$P_{\tau}(\mbx)_k > P_{\tau}(\mbx)_{\ell}$.

\bl                      \label{Muirhead:lemma:rearrange}
For all $n \in \N$ and  $\mbb   \in \Rn$, there exists a finite sequence of vectors 
$\mbb = \mbb_0, \mbb_1, \ldots, \mbb_r = \mbb^{\uparrow}$ 
and a finite sequence of transpositions $\tau_1,\tau_2,\ldots, \tau_r \in S_n$ 
such that $\mbb_i = P_{\tau_i}(\mbb_{i-1})$ and $\tau_i$ increases 
the vector $\mbb_{i-1}$ for all $i \in \{1,\ldots, r\}$.
\el

\begin{proof}
The proof is by induction on $n$.
If $n = 1$ and $\mbb = \bmat b_1 \emat$ 
or if $n = 2$ and $\mbb = \bsmallmat b_1 \\ b_2 \esmallmat$ 
with $b_1 \leq b_2$, then $\mbb = \mbb^{\uparrow}$ and  there is nothing to prove.  
If  $\mbb = \bsmallmat b_1 \\ b_2 \esmallmat$ with
$b_1 > b_2$, then the transposition $\tau_1 = (1,2) \in S_2$ increases \mbb\ and 
$P_{\tau_1}(\mbb) = \bsmallmat b_2 \\ b_1 \esmallmat = \mbb^{\uparrow}$.

Let $n \geq 3$ and assume that the Lemma holds for vectors in $\R^{n-1}$. 
Let $\mbb = \vectorsmallbn \in \Rn$.  Choose   $ i_0 \in \{1,\ldots,n \}$ 
such that $b_{i_0} = \min\{ b_1, b_2,\ldots, b_n \}$. 
If $b_1 = b_{i_0}$, then apply the Lemma to the vector 
$\bsmallmat b_2 \\ b_3 \\ \vdots \\ b_n \esmallmat \in \R^{n-1}$. 

If $b_1 > b_{i_0}$, then the transposition $\tau_1 = (1,i_0)$ increases \mbb. 
Let $\mbb_1 = P_{\tau_1}(\mbb)$. We have $ b_{\tau(1)} = b_{i_0}  \leq  b_{\tau(i)} $ for all 
$i \in \{1,\ldots, n\}$. Applying the Lemma to the vector 
$\bsmallmat b_{\tau(2)} \\ \vdots \\ b_{\tau(n)} \esmallmat \in \R^{n-1}$ completes the proof.  
\end{proof}

The \emph{standard inner product} of vectors 
$\mbx = \left( \begin{smallmatrix} x_1 \\ \vdots \\ x_n \end{smallmatrix}\right)$ 
and  $\mby = \vectorsmallyn$ is 
\[
\left( \mbx,\mby \right) = \sum_{j=1}^n x_j y_j.
\]
For all $\rho,\sigma \in S_n$ we have 
\begin{align*}
\left( P_{\rho}(\mbx),P_{\sigma}(\mby)\right) 
& = \sum_{j=1}^n x_{\rho^{-1}(j)} y_{\sigma^{-1}(j)} \\
& = \sum_{j=1}^n x_j y_{\sigma^{-1}\rho(j)} \\
& = \left( \mbx,P_{\rho^{-1}\sigma}(\mby)\right) 
\end{align*}
and so 
\[
\left( P_{\sigma}(\mbx),P_{\sigma}(\mby)\right) = \left( \mbx,\mby \right). 
\]

\bl                                  \label{Muirhead:lemma:rearrange-2}
Let \mba\ and \mbb\ be vectors in  \Rn\ and let $\tau$ be a transposition in $S_n$.
If  $\mba $ is increasing and $\tau$ increases \mbb, 
or if  $\mba $ is decreasing and  $\tau $ decreases \mbb, then 
\[
\left( \mba, P_{\tau}(\mbb) \right)  \geq  \left( \mba,\mbb \right). 
\]
If  $\mba $ is increasing and $\tau$ decreases \mbb, 
or if  $\mba $ is decreasing and  $\tau $ increases \mbb, then 
\[
\left( \mba, P_{\tau}(\mbb) \right)  \leq  \left( \mba,\mbb \right). 
\]
\el

\begin{proof}
Let $\mba= \vectorsmallan$ and $\mbb = \vectorsmallbn$, 
and let $\tau = (k,\ell)$ with $k < \ell$.     
Because $\tau^{-1} = \tau$ and $b_{\tau(i)} = b_i$ for all $i \neq k, \ell$, 
we have 
\begin{align*}
\left( \mba, P_{\tau}(\mbb) \right) - \left( \mba,\mbb \right) 
& = \sum_{i=1}^n a_i b_{\tau(i)} -  \sum_{i=1}^n a_i b_i \\
& = \sum_{i=1}^n a_i \left( b_{\tau(i)} -  b_i \right) \\
& = a_k (b_{\ell} - b_k) + a_{\ell}( b_k - b_{\ell}) \\
 & =   \left(  a_{\ell}  - a_k \right)  \left( b_k - b_{\ell} \right).  
\end{align*}
If \mba\ is increasing and $\tau$ increases \mbb, then $a_{\ell} \geq a_k$ 
and $b_k > b_{\ell}$, and so $\left( \mba, P_{\tau}(\mbb) \right) - \left( \mba,\mbb \right) \geq 0$. 
If \mba\ is decreasing and $\tau$ decreases \mbb, then $a_{\ell} \leq a_k$ 
and $b_k < b_{\ell}$, and so $\left( \mba, P_{\tau}(\mbb) \right) - \left( \mba,\mbb \right) \leq 0$. 
The proof of the second inequality is similar.  
\end{proof}

Let \mba\ and \mbb\ be vectors in \Rn.
Consider the set of all inner products of rearrangements of the vectors \mba\ and \mbb, 
that is, the set 
\[
S = \left\{ \left( P_{\sigma}(\mba), P_{\tau}(\mbb) \right): \sigma, \tau \in S_n \right\}.
\]
We have 
\[
S =  \left\{ \left( \mba, P_{\tau}(\mbb) \right):  \tau \in S_n \right\} 
=  \left\{ \left( \mba^{\uparrow}, P_{\tau}(\mbb) \right):  \tau \in S_n \right\}.
\]
The following is an inequality of Hardy, Littlewood, and P\' olya.

\bt
For all  vectors $\mba, \mbb \in \R^n$, 
\[
\left( \mba^{\uparrow}, \mbb^{\downarrow} \right) \leq  \left( \mba , \mbb  \right) 
\leq  \left( \mba^{\uparrow}, \mbb^{\uparrow} \right) 
= \left( \mba^{\downarrow}, \mbb^{\downarrow} \right) .
\]
\et

\begin{proof}
Let $\mba = \vectorsmallan$ and $\mbb = \vectorsmallbn$. 
Choose $\rho, \sigma \in S_n$ such that 
\[
\mba^{\uparrow} = P_{\rho}(\mba) = \bmat a_{\rho^{-1}(1)}\\ a_{\rho^{-1}(2)}\\ \vdots \\a_{\rho^{-1}(n)}\emat
\qqand 
\mbb^{\uparrow} = P_{\sigma}(\mbb) = \bmat b_{\sigma^{-1}(1)}\\ b_{\sigma^{-1}(2)}\\ \vdots \\b_{\sigma^{-1}(n)}\emat. 
\]
Then
\[
\mba^{\downarrow} =  \bmat a_{\rho^{-1}(n)} \\ \vdots \\a_{\rho^{-1}(2)}\\ a_{\rho^{-1}(1)} \emat
\qqand 
\mbb^{\downarrow} =  \bmat a_{\rho^{-1}(n)} \\ \vdots \\ b_{\sigma^{-1}(2)}\\ b_{\sigma^{-1}(1)} \emat 
\]
and so
\[
\left( \mba^{\uparrow}, \mbb^{\uparrow} \right)   
 = \sum_{j=1}^n a_{\rho^{-1}(j)}a_{\rho^{-1}(j)}
 =  \left( \mba^{\downarrow}, \mbb^{\downarrow} \right).
\]

Next we prove that 
\[
 \left( \mba^{\uparrow}, \mbb  \right) 
\leq  \left( \mba^{\uparrow}, \mbb^{\uparrow} \right). 
\]
By Lemma~\ref{Muirhead:lemma:rearrange}, there is a finite sequence of vectors 
$\mbb = \mbb_0, \mbb_1, \ldots, \mbb_r = \mbb^{\uparrow}$ 
and a finite sequence of transpositions $\tau_1,\tau_2,\ldots, \tau_r \in S_n$ 
such that $\mbb_i = P_{\tau_i}(\mbb_{i-1})$ and $\tau_i$ increases 
 $\mbb_{i-1}$ for all $i \in \{1,\ldots, r\}$.
By Lemma~\ref{Muirhead:lemma:rearrange-2}, 
\begin{align*}
\left( \mba^{\uparrow}, \mbb^{\uparrow} \right) -  \left( \mba^{\uparrow} , \mbb  \right) 
& = \left( \mba^{\uparrow}, \mbb_r \right) -  \left( \mba^{\uparrow}, \mbb_0  \right) \\
& = \sum_{i=1}^r \left( \left( \mba^{\uparrow}, \mbb_i \right) - \left( \mba^{\uparrow}, \mbb_{i-1} \right) \right) \\
& = \sum_{i=1}^r \left( \left( \mba^{\uparrow}, P_{\tau_i}(\mbb_{i-1}) \right) - \left( \mba^{\uparrow},\mbb_{i-1}  \right) \right) \\
& \geq 0.
\end{align*}
It follows that 
\[
\left( \mba , \mbb  \right)  = \left( P_{\rho}(\mba), P_{\rho}(\mbb) \right) 
=  \left( \mba^{\uparrow}, P_{\rho}(\mbb)  \right) 
\leq  \left( \mba^{\uparrow}, P_{\rho}(\mbb)^{\uparrow}  \right) 
= \left( \mba^{\uparrow}, \mbb^{\uparrow}  \right).  
\]  

There is a similar  proof of the inequality 
\[
 \left( \mba^{\uparrow}, \mbb^{\downarrow}  \right) 
\leq  \left( \mba^{\uparrow}, \mbb \right).   
\]
By Lemma~\ref{Muirhead:lemma:rearrange}, there is a finite sequence of vectors 
$\mbb = \mbb_0, \mbb_1, \ldots, \mbb_r = \mbb^{\uparrow}$ 
and a finite sequence of transpositions $\tau_1,\tau_2,\ldots, \tau_r \in S_n$ 
such that $\mbb_i = P_{\tau_i}(\mbb_{i-1})$ and $\tau_i$ decreases 
 $\mbb_{i-1}$ for all $i \in \{1,\ldots, r\}$.
By Lemma~\ref{Muirhead:lemma:rearrange-2}, 
\begin{align*}
\left( \mba^{\uparrow}, \mbb^{\uparrow} \right) -  \left( \mba^{\uparrow} , \mbb  \right) 
& = \left( \mba^{\uparrow}, \mbb_r \right) -  \left( \mba^{\uparrow}, \mbb_0  \right) \\
& = \sum_{i=1}^r \left( \left( \mba^{\uparrow}, \mbb_i \right) - \left( \mba^{\uparrow}, \mbb_{i-1} \right) \right) \\
& = \sum_{i=1}^r \left( \left( \mba^{\uparrow}, P_{\tau_i}(\mbb_{i-1}) \right) - \left( \mba^{\uparrow},\mbb_{i-1}  \right) \right) \\
& \leq 0.
\end{align*}
It follows that 
\[
 \left( \mba^{\uparrow}, \mbb^{\downarrow}  \right) 
 = \left( \mba^{\uparrow}, P_{\rho}(\mbb)^{\downarrow} \right) 
\leq  \left( \mba^{\uparrow}, P_{\rho}(\mbb)  \right) 
= \left( P_{\rho}(\mba) , P_{\rho}(\mbb)  \right) 
= \left( \mba, \mbb  \right).  
\]  
This completes the proof. 
\end{proof}

\section{Vector majorization and doubly stochastic matrices} 

Let $\R_{\geq 0}$ denote the set of nonnegative real numbers.  
A \emph{nonnegative vector} is a vector with nonnegative coordinates.  
The \emph{nonnegative octant}\index{octant} in \Rn\ is the set 
of nonnegative vectors in $\R^n$, that is, 
\[
\Rn_{\geq 0} = \left\{ \mba = \vectoran \in \R^n: a_i \geq 0 \text{ for all } i \in \{ 1,\ldots, n \} \right\}. 
\]
If $\mba \in \Rn_{\geq 0}$, then  $P_{\sigma}(\mba) \in \Rn_{\geq 0}$ for all $\sigma \in S_n$. 
In particular, if $\mba \in \Rn_{\geq 0}$, then $\mba^{\downarrow} \in \Rn_{\geq 0}$ 
and $\mba^{\uparrow} \in \Rn_{\geq 0}$. 

Let $\mba$ and $\mbb$ be vectors in $\Rn$ and let 
\[
\mba^{\downarrow}  = \bmat a_1^{\downarrow} \\ \vdots \\ a_n^{\downarrow} \emat 
\qqand \mbb^{\downarrow}  = \bmat b_1^{\downarrow} \\ \vdots \\ b_n^{\downarrow} \emat
\] 
be the associated decreasing vectors.  
Thus, 
\[
a_1^{\downarrow} \geq \cdots \geq a_n^{\downarrow} 
\qqand 
b_1^{\downarrow} \geq \cdots \geq b_n^{\downarrow}. 
\]
The  vector $\mba$ \emph{majorizes}\index{majorize} 
the  vector $\mbb$  
if the vectors $\mba$ and $\mbb$ have the same coordinate sum 
and if 
\[
\sum_{i=1}^k b_i^{\downarrow} \leq \sum_{i=1}^k a_i^{\downarrow} 
\qquad \text{for all $k \in \{1,\ldots, n-1\}$}.  
\]
We write 
\[
\mbb \prec \mba
\]
if the vector \mbb\ is majorized by the vector  \mba. 

Majorization is a partial order  but  not a total order on 
the set of vectors in $\Rn_{\geq 0}$ with the same coordinate sum.

For example, in $\R^3_{\geq 0}$, we have 
\[
\bmat 5 \\ 4 \\ 3 \emat \prec \bmat 7 \\ 3 \\ 2 \emat \qqand 
\bmat 2 \\ 6  \\ 5  \emat \prec \bmat  4 \\ 1 \\ 8  \emat.
\]
The vectors $\mba = \bsmallmat 7 \\ 5\\ 1 \esmallmat$ 
and $\mbb = \bsmallmat 8 \\ 3 \\ 2 \esmallmat$ 
are decreasing with coordinate sum 13,  but neither vector majorizes the other.  

For every positive integer $n$ we have the following linearly ordered majorization chain: 
\[
\bmat 1/n \\ 1/n \\ 1/n  \\ \vdots \\ 1/n  \\ 1/n \emat \prec 
\bmat 1/(n-1) \\ 1/(n-1)  \\ 1/(n-1)\\ \vdots \\ 1/(n-1) \\ 0 \emat \prec \cdots 
\prec \bmat 1/3 \\ 1/3 \\ 1/3  \\ 0  \\ \vdots \\ 0 \emat 
\prec \bmat 1/2 \\ 1/2 \\ 0  \\ 0  \\ \vdots \\ 0 \emat \prec \bmat 1 \\ 0  \\ 0  \\ 0  \\ \vdots \\ 0 \emat. 
\]

We have the following majorization relations.

\bl       \label{Muirhead:lemma:Chong-1} 
If $i_1,\ldots, i_k$ are integers such that 
$1 \leq i_1 < i_2 < \cdots < i_k \leq n$ 
and if $x_1,\ldots, x_n$ be real numbers such that $x_1 \geq x_2  \geq  \cdots  \geq x_n$, 
then 
\[
\sum_{j=1}^k x_{i_j} \leq \sum_{j=1}^k x_j.
\]  
\el

\begin{proof}
We have $i_1 \geq 1$.  
If $j \in \{1,\ldots, k-1\}$ and $i_j \geq j$, then $i_{j+1} > i_j \geq j$ and so  $i_{j+1} \geq j+1$.
Therefore, $i_j \geq j $ for all $j \in \{1,\ldots, k\}$ and so
$x_{i_j} \leq x_j$  for all $j \in \{1,\ldots, k\}$. 
The inequality follows immediately.  
\end{proof}

The following majorization inequality is due to Chong~\cite{chon74}.

\bl                    \label{Muirhead:lemma:Chong-4}
For all integers $\ell \geq 2$ and all vectors $\mba_1, \mba_2, \ldots, \mba_{\ell} \in \Rn$, 
\[
\left( \mba_1 + \mba_2 + \cdots+  \mba_{\ell}  \right)^{\downarrow} \prec 
\mba_1^{\downarrow} + \mba_2^{\downarrow} + \cdots \mba_{\ell}^{\downarrow}.
\]
\el 

\begin{proof}
For all $j \in \{1,\ldots, \ell\}$, let 
\[
\mba_j = \bmat a_{1,j} \\ \cdots \\ a_{n,j} \emat
\]
and choose $\rho_i \in S_n$ such that 
\[
\mba_j^{\downarrow} = \bmat a_{ \rho^{-1}(1)},j \\ \cdots \\ a_{\rho^{-1}(n),j} \emat.
\]
Choose $\sigma \in S_n$ such that 
\[
\left( \mba_1 + \mba_2 + \cdots+  \mba_{\ell}  \right)^{\downarrow} 
= \bmat
\sum_{j=1}^{\ell}  a_{ \sigma^{-1}(1),j} \\ \vdots \\ \sum_{j=1}^{\ell}  a_{\sigma^{-1}(n),j}   \\
\emat.
\]
Let 
\[
\mbb=\left( \mba_1 + \mba_2 + \cdots+  \mba_{\ell}  \right)^{\downarrow} = \vectorbn
\]
and 
\[
\mbc = \mba_1^{\downarrow} + \mba_2^{\downarrow} + \cdots \mba_{\ell}^{\downarrow} = \vectorcn. 
\]

The coordinates in the vectors $\mba_1^{\downarrow}, \ldots, \mba_{\ell}^{\downarrow}$ 
are decreasing and so, by Lemma~\ref{Muirhead:lemma:Chong-1}, 
\[
\sum_{i=1}^k a_{\sigma^{-1}(i),j}  \leq \sum_{i=1}^k a_{\rho^{-1}(i),j}  
\]
for all $j \in \{1,\ldots, \ell\}$ and  $k \in \{1,\ldots, n-1\}$ 
and so 
\begin{align*}
\sum_{i=1}^k b_i 
& =  \sum_{i=1}^k  \sum_{j=1}^{\ell} a_{\sigma^{-1}(i),j}  
= \sum_{j=1}^{\ell}  \sum_{i=1}^k   a_{\sigma^{-1}(i),j} \\ 
& \leq  \sum_{j=1}^{\ell}  \sum_{i=1}^k  a_{\rho^{-1}(i),j} 
=   \sum_{i=1}^k  \sum_{j=1}^{\ell} a_{\rho^{-1}(i),j} \\
&  = \sum_{i=1}^k c_j.
\end{align*}
This completes the proof.  
\end{proof}

\bl                    \label{Muirhead:lemma:Chong-5}
Let $\mba_1, \mba_2, \ldots, \mba_{\ell} \in \Rn$ 
and $\mbb_1, \mbb_2, \ldots, \mbb_{\ell} \in \Rn$   
with $\mbb_i \prec \mba_i$ for all $i \in \{1,\ldots, \ell\}$. 
If $\lambda_1,\ldots, \lambda_{\ell} \in \R_{\geq 0}$, then 
\[
\lambda_1  \mbb_1 + \lambda_2 \mbb_2 + \cdots+ \lambda_{\ell}  \mbb_{\ell} 
\prec \lambda_1 \mba_1 + \lambda_2 \mba_2 
+ \cdots \lambda_{\ell}  \mba_{\ell}.
\]
Let $\mba, \mbb_1, \mbb_2, \ldots, \mbb_{\ell} \in \Rn$   
with $\mbb_i \prec \mba$ for all $i \in \{1,\ldots, \ell\}$. 
If $\lambda_1,\ldots, \lambda_{\ell} \in \R_{\geq 0}$ 
and $\sum_{i=1}^{\ell} \lambda_i = 1$,  then 
\[
\lambda_1  \mbb_1 + \lambda_2 \mbb_2 + \cdots+ \lambda_{\ell}  \mbb_{\ell} 
\prec   \mba.
\]
\el 

\begin{proof}
This follows from Lemma~\ref{Muirhead:lemma:Chong-4}.
\end{proof}

In \Rn, for all $j \in \{1,\ldots, n\}$ we have 

The standard basis vectors $\{\mbe_1,\ldots, \mbe_n\}$ in $\R^n$ are defined by   
$\mbe_j = \bsmallmat \delta_{1,j} \\ \vdots \\ \delta_{n,j} \esmallmat$ for all $j \in \{1,\ldots, n\}$, 
where $\delta_{i,j}$ is the Kronecker delta.   
Let $\mbf_n = \bsmallmat 1/n  \\ \vdots \\ 1/n \esmallmat \in \Rn$. 
We have $\mbe_j^{\downarrow} = \mbe_1$ for all $j \in \{1,\ldots, n\}$ 
and $\mbf_n^{\downarrow} = \mbf_n$.

\bt                  \label{Muirhead:theorem:jn-a-e} 
For every vector $\mba = \vectorsmallan$  in $\Rn_{\geq 0} $ with $\sum_{i=1}^n a_i = 1$, 
there is the majorization inequality   
\[
\mbf_n \prec \mba \prec \mbe_k 
\]
for all $k \in \{1,\ldots, n\}$.  Moreover, 
\[
\mba \prec \mbf_n \quad \text{implies} \quad \mba = \mbf_n 
\]
and 
\[
\mbe_1 \prec \mba \quad \text{implies $\mba = \mbe_k$ 
for some $k \in \{1,\ldots, n\}$.}
\]
\et

\begin{proof}
Let $\mba^{\downarrow} = \bsmallmat a_1^{\downarrow}  \\ \vdots \\ a_n^{\downarrow} \esmallmat$, 
where $a_1^{\downarrow}  \geq \cdots \geq a_n^{\downarrow} \geq 0$ 
and $ \sum_{i=1}^n a_i^{\downarrow} =1$. 
For all $k \in \{1,2,\ldots, n\}$,  we have 
$\mbe_k^{\downarrow} =\mbe_1 = \bsmallmat \delta_{1,1} \\ \vdots \\ \delta_{n,1} \esmallmat$ 
and 
\[
\sum_{i=1}^k a_i^{\downarrow} \leq \sum_{i=1}^n a_i^{\downarrow} = 1 = \delta_{1,1} = \sum_{i=1}^k \delta_{i,1} 
\]
and so $\mba \prec \mbe_1 = \mbe_k^{\downarrow}$.   

If 
\[
\sum_{i=1}^{k} a_i^{\downarrow} < \frac{k}{n} 
\]
for some $k \in \{1,\ldots, n\}$, then 
\[
ka_k^{\downarrow}  \leq \sum_{i=1}^{k} a_i^{\downarrow} < \frac{k}{n} 
\]
implies 
\[
a_k^{\downarrow}  < \frac{1}{n} 
\]
and so 
\[
1 = \sum_{i=1}^n a_i^{\downarrow} = \sum_{i=1}^k a_i^{\downarrow} + \sum_{i=k+1}^n a_i^{\downarrow} 
 < \frac{k}{n} +  \frac{n-k}{n} = 1
\]
which is absurd.  Therefore, $k/n \leq  \sum_{i=1}^k a_i^{\downarrow}$ 
for all $k \in \{1,\ldots, n\}$ and   $\mbf_n \prec \mba$.

If $\mba \prec \mbf_n$, then,   for all $k \in \{1,\ldots, n\}$, we have 
\[
ka_k^{\downarrow} \leq \sum_{i=1}^k a_i^{\downarrow} \leq \frac{k}{n}
\]
and so
\[
a_k^{\downarrow} \leq \frac{1}{n}.   
\]
The coordinate sum $\sum_{k=1}^n a_k^{\downarrow} = 1$ implies that 
$a_k^{\downarrow} = 1/n$ for all $k$ and  
$\mba = \mbf_n$.   

If $\mbe_1 \prec \mba$, 
then $1 = \delta_{1,1} \leq a^{\downarrow}_1 \leq 1$ and so  $a^{\downarrow}_1 = 1$.   
It follows from $\sum_{i=1}^n a_i^{\downarrow} = 1$ that 
$a^{\downarrow}_i = 0$ for all $i \in \{2,3,\ldots, n\}$ and so  
$\mba^{\downarrow} = \mbe_1$ and $\mba =  \mbe_k$ 
for some $k \in \{1,\ldots, n\}$.  
This completes the proof. 
\end{proof}

\bl                 \label{Muirhead:lemma:DS-coordinateSum} 
Let $P = \bmat p_{i,j} \emat$ be a column stochastic $n\times n$ matrix.  
If $\mba = \vectorsmallan \in \Rn_{\geq 0}$ and $P(\mba) = \mbb = \vectorsmallbn$, then
$\mbb  \in \Rn_{\geq 0}$ and 
\[
 \sum_{i=1}^n b_i = \sum_{j=1}^n a_j. 
\]
\el

\begin{proof}
For all $i \in \{1,\ldots, n\}$, we have 
\[
b_i  = \sum_{j=1}^n p_{i,j} a_j \geq 0 
\]
and so $\mbb  \in \Rn_{\geq 0}$.  Moreover, 
\[
\sum_{i=1}^n b_i = \sum_{i=1}^n \sum_{j=1}^n p_{i,j} a_j 
= \sum_{j=1}^n a_j  \sum_{i=1}^n p_{i,j} = \sum_{j=1}^na_j  \colsum_j(P) = \sum_{j=1}^n a_j. 
\]
This completes the proof. 
\end{proof}

\bt                \label{Muirhead:theorem:DoublyStochasticMajorize}
Let $P = \bmat p_{i,j} \emat$ be a nonnegative $n\times n$ matrix. 
The following are equivalent:
\benum
\item[(i)]
$P(\mba) \prec \mba$ for all vectors $\mba \in \Rn_{\geq 0}$. 
\item[(ii)]
$P(\mbe_j) \prec \mbe_j$ for all   $j \in \{1,\ldots, n \}$ and $P(\mbf_n) \prec \mbf_n$ . 
\item[(iii)]
The matrix $P = \bmat p_{i,j} \emat$ is doubly stochastic. 
\eenum
\et

\begin{proof}
Clearly, (i) implies~(ii).

Assume~(ii).  For all  $j \in \{1,\ldots, n\}$, the $j$th column of the matrix $P$ is 
\[
\col_j(P) = \bmat p_{1,j} \\ \vdots \\ p_{n,j} \emat = P(\mbe_j ).  
\]
If $P(\mbe_j ) \prec \mbe_j $, then 
\[
1 = \colsum(\mbe_j) = \colsum_j(P) = \sum_{i=1}^n p_{i,j} 
\]
and so  $P$ is column stochastic.

By Theorem~\ref{Muirhead:theorem:jn-a-e}, the majorization relation 
$P(\mbf_n) \prec \mbf_n$ implies $P(\mbf_n) = \mbf_n$.  
The $i$th coordinate of the vector $\mbf_n$ is $1/n$ 
and the $i$th coordinate of the vector $\mbf_n$ is $\sum_{j=1}^n \frac{p_{i,j}}{n} $.  
Therefore, 
\[
 \frac{1}{n} =   \sum_{j=1}^n \frac{p_{i,j}}{n} =  \frac{1}{n} \rowsum_i(P) 
\]
and so $\rowsum_i(P) = 1$ for all $i \in \{1,\ldots, n\}$  
and the matrix $P$ is also row stochastic.  
Thus, (ii) implies~(iii).

Let $P$ be a doubly stochastic matrix.  
 Let $\mba \in \R^n_{\geq 0}$ and $\mbb = P(\mba)$. 
 We shall prove that $\mbb^{\downarrow}  \prec \mba^{\downarrow}$.

There exist permutations $\sigma$ and $\tau$ in $S_n$ such that 
\[
P_{\sigma}(\mba) = \mba^{\downarrow} \qqand P_{\tau}(\mbb) = \mbb^{\downarrow}
\]
and so 
\[
\mbb^{\downarrow} = P_{\tau}(\mbb) = P_{\tau}P(\mba) 
= P_{\tau}PP_{\sigma^{-1}}\left( \mba^{\downarrow}\right) = Q\mba^{\downarrow}.  
\]
The matrices $P_{\tau}$ and $P_{\sigma^{-1}}$ are doubly stochastic 
and so the matrix $Q = P_{\tau}PP_{\sigma^{-1}} = \bmat q_{i,j} \emat$ is doubly stochastic.

For $k \in \{1,\ldots, n-1\}$ and $j \in \{1,\ldots, n \}$, define the partial column sums 
\[
c_j^{(k)} =  \sum_{i=1}^k  q_{i,j}. 
\]
We have 
\[
0 \leq c_j^{(k)}  \leq  \sum_{i=1}^n  q_{i,j} = \colsum_j(Q) = 1 
\]
and 
\[
 \sum_{j=1}^n  c_j^{(k)} =  \sum_{j=1}^n \sum_{i=1}^k  q_{i,j} 
 = \sum_{i=1}^k  \sum_{j=1}^n  q_{i,j}  = \sum_{i=1}^k \rowsum_i(Q) = k.
\]
Let
\[
\mba^{\downarrow} = \bmat a_1^{\downarrow} \\ \vdots \\ a_n^{\downarrow} \emat 
\qqand \mbb^{\downarrow} = \bmat b_1^{\downarrow} \\ \vdots \\ b_n^{\downarrow} \emat. 
\]
From Lemma~\ref{Muirhead:lemma:DS-coordinateSum} we obtain 
$\sum_{i=1}^n b_i^{\downarrow} = \sum_{i=1}^n a_i^{\downarrow}$.   
The inequalities 
$a_1^{\downarrow} \geq a_2^{\downarrow} \geq \cdots \geq a_n^{\downarrow}$ imply that  
\[
b_1^{\downarrow} = \sum_{j=1}^n q_{1,j} a_j^{\downarrow} 
 \leq a_1^{\downarrow}  \sum_{j=1}^n q_{1,j} 
 = a_1^{\downarrow} \rowsum_1(Q) = a_1^{\downarrow}.  
\]
For all $k \in \{2,\ldots, n-1\}$ we have 
\begin{align*}
\sum_{i=1}^k b_i^{\downarrow} & = \sum_{i=1}^k \sum_{j=1}^n q_{i,j}a_j^{\downarrow} 
= \sum_{j=1}^n  a_j^{\downarrow}  \sum_{i=1}^k  q_{i,j} =  \sum_{j=1}^n a_j^{\downarrow} c_j^{(k)}  \\
& \leq  \sum_{j=1}^{k-1} c_j^{(k)}  a_j^{\downarrow} + a_k^{\downarrow} \sum_{j= k}^n c_j^{(k)}  
=  \sum_{j=1}^{k-1} c_j^{(k)} a_j^{\downarrow} + a_k^{\downarrow} \left( k - \sum_{j= 1}^{k-1} c_j^{(k)} \right) \\
&=  \sum_{j=1}^{k-1} c_j^{(k)} (a_j^{\downarrow} - a_k^{\downarrow})  + k a_k^{\downarrow}  
 \leq  \sum_{j=1}^{k-1} (a_j^{\downarrow} - a_k^{\downarrow})  + k a_k^{\downarrow} \\
& = \sum_{j=1}^k a_j^{\downarrow}. 
\end{align*}
Therefore,  $\mbb \prec \mba$ and (iii) implies (i).
This completes the proof. 
\end{proof}

\section{The Hardy-Littlewood-P\' olya $T$-transformation} 

Let $\mba, \mbb \in \Rn$.  
In  Theorem~\ref{Muirhead:theorem:DoublyStochasticMajorize} 
we proved that if $P$ is an $n \times n$ doubly stochastic  matrix, then 
\mba\ majorizes $P(\mba)$, that is,  if $\mbb = P(\mba)$, then 
$\mbb \prec \mba$.  
In this section we prove the converse: 
If $\mbb \prec \mba$, 
then there is a doubly stochastic matrix $P$ such that $\mbb = P(\mba)$.  
The proof is based on the following elementary lemma.   

\bl                             \label{Muirhead:lemma:HardyLittlewoodPolya-T0}
Let $a_k, a_{\ell}, c_k$, and $c_{\ell}$ be real numbers such that 
\beq                                \label{Muirhead:HLP-1}
a_{\ell} + a_k =  c_{\ell} + c_k
\eeq
and
\beq                                \label{Muirhead:HLP-2}
a_{\ell} < c_{\ell} \leq c_k < a_k.
\eeq
There is a unique $2 \times 2$ doubly stochastic matrix $T_{k,\ell}$ such that 
\[
T_{k,\ell}\bmat a_k\\ a_{\ell} \emat = \bmat c_k\\ c_{\ell} \emat.
\] 
Moreover, 
\[
T_{k,\ell} = \bmat\lambda  & 1 - \lambda \\ 1 - \lambda & \lambda \emat 
\]
where 
\[
\lambda  = \frac{c_k-a_{\ell}}{a_k-a_{\ell}} = \frac{a_k-c_{\ell}}{a_k-a_{\ell}}.
\]
\el

\begin{proof}
Every $2 \times 2$ doubly stochastic matrix is of the form 
$\bmat\lambda  & 1 - \lambda \\ 1 - \lambda & \lambda \emat$ 
for some $\lambda \in [0,1]$.  We have 
\[
\bmat\lambda  & 1 - \lambda \\ 1 - \lambda & \lambda \emat 
 \bmat a_k \\ a_{\ell} \emat 
= \bmat c_k \\ c_{\ell} \emat
\]
if and only if 
\[
\lambda a_k + (1-\lambda) a_{\ell} = c_k 
\]
and 
\[
(1-\lambda) a_k + \lambda a_{\ell} = c_{\ell}. 
\]
The unique solution of the first equation is 
\[
\lambda_1 = \frac{c_k-a_{\ell}}{a_k-a_{\ell}} 
\]
and the unique solution of the second equation is 
\[
\lambda_2 = \frac{a_k-c_{\ell}}{a_k-a_{\ell}}.
\]
Equation~\eqref{Muirhead:HLP-1} implies that 
$\lambda_1 = \lambda_2 = \lambda$.   
Inequality~\eqref{Muirhead:HLP-2} implies that $0 < \lambda < 1$. 
This completes the proof.  
\end{proof}

For example, if $\bmat a_2\\ a_5 \emat = \bmat 4 \\ -10 \emat$ and 
$\bmat c_2\\ c_5 \emat = \bmat 1 \\ -7\emat$, then $\lambda = 11/14$ and 
\[
T_{2,5} \bmat 4 \\ -10 \emat = 
\bmat 11/14 & 3/14 \\ 3/14 & 11/14 \emat  \bmat 4 \\ -10 \emat = \bmat 1 \\ -7\emat. 
\]

\bl                               \label{Muirhead:lemma:HardyLittlewoodPolya-T1}
Let $k < \ell$ and let $\mba = \vectorsmallan$ and  $\mbc = \vectorsmallcn$ be  vectors in \Rn\ such that 
 \[
a_{\ell} + a_k =  c_{\ell} + c_k
\]
\[
a_{\ell} < c_{\ell} \leq c_k < a_k
\]
and 
\[
c_i = a_i \text{ for all $i \neq k,\ell$.}
\]
There is a doubly stochastic $n \times n$ matrix $T = \bmat t_{i,j} \emat$ 
such that $T\mba = \mbc$.  
\el

\begin{proof}
Let 
\[
\lambda = \frac{c_k-a_{\ell}}{a_k-a_{\ell}}  = \frac{a_k-c_{\ell}}{a_k-a_{\ell}}.
\]
Let $T = \bmat t_{i,j}\emat $ be the $n \times n$ matrix defined as follows: 
\[
t_{k,k} = t_{\ell,\ell} =  \lambda
\]
\[
t_{k,\ell} = t_{\ell,k} =  1 - \lambda 
\]
and
\[
t_{i,j} = \delta_{i,j} \qquad \text{for all $(i,j) \notin  \{ (k,k), (k,\ell), (\ell,k), (\ell,\ell) \}$}
\]
where $\delta_{i,j}$ is the Kronecker delta.  
The matrix $T$ looks like 
\[
\setcounter{MaxMatrixCols}{14} 
T = 
\bmat 
1 & \cdots &  0 & 0 & 0 & \cdots &  0 & 0 & 0  & \cdots & 0  \\
\vdots &&&&&&&&& & \vdots  \\
0 & \cdots &  1 & 0 & 0 & \cdots &  0 & 0 & 0  & \cdots  & 0 \\
 0 & \cdots &  0 & \lambda & 0 &  \cdots &  0 & 1 - \lambda & 0 & \cdots & 0   \\
0 & \cdots &  0 & 0 & 1 & \cdots &  0 & 0 & 0  & \cdots & 0  \\
\vdots &&&&&&&&&& \vdots  \\
0 & \cdots &  0 & 0 & 0 & \cdots &  1 & 0 & 0  & \cdots  & 0 \\
 0 & \cdots &  0 & 1 - \lambda & 0 & \cdots &0 & \lambda & 0 & \cdots & 0  \\
 0 & \cdots &  0 & 0 & 0 & \cdots &  0 & 0 & 1  & \cdots & 0  \\ 
 \vdots &&&&&&&&&& \vdots  \\
  0 & \cdots &  0 & 0 & 0 & \cdots &  0 & 0 & 0  & \cdots & 1  \\
\emat. 
\]
It follows from Lemma~\ref{Muirhead:lemma:HardyLittlewoodPolya-T0} 
that $T\mba = \mbc$.  
One sees directly that $T$ is doubly stochastic.  
This completes the proof.  
\end{proof}

The matrix $T$ constructed in Lemma~\ref{Muirhead:lemma:HardyLittlewoodPolya-T1} 
is the \index{T-transformation}\index{Hardy-Littlewood-P\' olya $T$-transformation}
\emph{Hardy-Littlewood-P\' olya $T$-transformation}.

For example, if $n=5$, $k=2$, $\ell = 3$, and 
\[
\mba = \bmat 3 \\ 4 \\ 5 \\ 8\\ -10 \emat \qqand \mbb = \bmat 3 \\ 1 \\ 5 \\ 8\\ -7 \emat
\]
then
\[
T = \bmat 
1 & 0 & 0 & 0 & 0 \\
0 & 11/14 & 0 & 0 & 3/14 \\
0& 0 & 1 & 0 & 0 \\
0 & 0 & 0 & 1 & 0 \\
0 & 3/14 & 0 & 0 & 11/14 
\emat
\]
and $T\mba = \mbb$.

The \emph{Hamming distance} between vectors $\mba = \vectorsmallan$ 
and $\mbb =\vectorsmallbn$ in $\R^n$ is
\beq        \label{Muirhead:Hamming}
d_H(\mba,\mbb) = \card\{i \in \{ 1,\ldots, n\}: a_i  \neq b_i \}.
\eeq

\bl             \label{Muirhead:lemma:Hamming} 
For all vectors $\mba, \mbb, \mbc \in \Rn$, 
\benum
\item[(i)]
$d_H(\mba,\mbb) \in \{0,1,2,\ldots, n\}$ 
\item[(ii)]
$d_H(\mba,\mbb) = 0$ if and only if $\mba = \mbb$  
\item[(iii)]
Symmetry: 
$d_H(\mba,\mbb) =  d_H(\mbb,\mba)$
\item[(iv)]
Triangle inequality:
$d_H(\mba,\mbb) \leq d_H(\mba,\mbc) + d_H(\mbc,\mbb) $
\item[(v)]
If $\mba \neq \mbb$ and $\sum_{i=1}^n a_i = \sum_{i=1}^n b_i$, 
then $d_H(\mba,\mbb) \geq 2$. 
\eenum
\el

\begin{proof}
Properties (i), (ii), and (iii) are immediate.  

To prove (iv), which is the triangle inequality for vectors $\mba, \mbb ,\mbc$, let 
\begin{align*}
I & = \{i \in \{ 1,\ldots, n\}: a_i  \neq c_i \} \\
J & = \{j \in \{ 1,\ldots, n\}: c_j  \neq b_j \} \\
K & =   \{k \in \{ 1,\ldots, n\}: a_k  \neq b_k \} \\
L & =  \{1,\ldots, n\} \setminus (I \cup J).
\end{align*} 
If $\ell \in L$, then $\ell \notin I$ and $\ell \notin J$, and so $a_{\ell} = c_{\ell}$ and $c_{\ell} = b_{\ell}$,
hence $a_{\ell} = b_{\ell}$. 
If $k \in K$, then $a_k \neq b_k$ and so $k \notin L$, that is, 
$k \in  \{1,\ldots, n\} \setminus L = I \cup J$. 
 It follows that  $K \subseteq I \cup J$ and  
\[
d_H(\mba,\mbb) = |K| \leq |I\cup J| \leq |I| + |J| = d_H(\mba,\mbc) + d_H(\mbc,\mbb).
\]

To prove (v), we observe that if $\mba \neq \mbb$, then $a_{i_1} \neq b_{i_1}$ for some 
$i_1 \in \{1,2,\ldots, n\}$.  If $\sum_{i=1}^n a_i = \sum_{i=1}^n b_i$, then 
\[
\sum_{\substack{i=1\\i \neq i_1}}^n a_i \neq \sum_{\substack{i=1\\i \neq i_1}}^n b_i 
\]
and so $a_{i_2} \neq b_{i_2}$ for some 
$i_2 \in \{1,2,\ldots, n\} \setminus \{i_1\}$.  Therefore, $d_H(\mba,\mbb) \geq 2$.  
This completes the proof.  
\end{proof}

\bl                                     \label{Muirhead:lemma:HardyLittlewoodPolya-Tc}
Let \mba\ and \mbb\ be  decreasing  vectors in $\R^n$  such that $\mbb \prec \mba$ 
and $d_H(\mba,\mbb) > 2$. 
There exists a vector \mbc\ in $\R^n$  such that 
\benum
\item[(i)]
\mbc\ is decreasing,
\item[(ii)]
$\mbb \prec \mbc \prec \mba$, 
\item[(iii)]
$d_H(\mba,\mbc) = 2$ and $2 \leq d_H(\mbb,\mbc) < d_H(\mbb,\mba)$, 
\item[(iv)]
$T\mba = \mbc$ for some doubly stochastic matrix $T$. 
\eenum
\el

\begin{proof}
Let  $\mba = \vectorsmallan$ and $\mbb = \vectorsmallbn$.  
We have $\mba  \neq \mbb $ because $d_H(\mba,\mbb) > 2$ 
and so there is a smallest integer  
$k_1$ such that $b_{k_1} \neq a_{k_1}$.
Therefore, $b_i = a_i$ for all $i \in \{1,\ldots, k_1-1\}$ and  
$\sum_{i=1}^{k_1-1} b_i = \sum_{i=1}^{k_1-1} a_i$.   
The inequality $\sum_{i=1}^{k_1} b_i \leq \sum_{i=1}^{k_1} a_i$ implies that 
$b_{k_1} \leq a_{k_1}$ and so $b_{k_1} < a_{k_1}$.
The equality $\sum_{i=1}^n b_i = \sum_{i=1}^n  a_i$ implies that $b_{\ell_1} > a_{\ell_1}$ 
for some integer $\ell_1  > k_1$.  

Let $\ell$ be the smallest integer such that $b_{\ell} > a_{\ell}$, and let $k$ be the largest integer  
such that $k < \ell$ and  $b_k < a_k$.  
Because the vectors \mba\ and \mbb\ are decreasing, we have 
\[
a_{\ell} < b_{\ell} \leq  b_k < a_k
\]
\[
a_i = b_i  \qquad \text{for all $i \in \{ k+1,k+2, \ldots, \ell-1\}$} 
\]
and 
\[
\delta = \min(a_k-b_k, b_{\ell} - a_{\ell} ) > 0. 
\]
Define the vector $\mbc = \vectorsmallcn$ as follows:
\begin{align*} 
c_k    & = a_k  - \delta \\
c_{\ell} & = a_{\ell}+\delta \\
c_i & = a_i \quad \text{for all $i \neq k,\ell$.} 
\end{align*} 
It follows that 
\[
c_1 \geq c_2 \geq \cdots \geq c_{k-1}
\]
\[
c_{k+1} \geq c_{k+2} \geq \cdots \geq c_{\ell-1}
\]
and 
\[
c_{\ell+1} \geq c_{\ell+2} \geq \cdots \geq c_n.
\]
We have 
\[
c_{k+1} = a_{k+1} = b_{k+1} \leq b_k \leq a_k - \delta = c_k < a_k \leq a_{k-1} = c_{k-1} 
\]
and
\[
c_{\ell+1} = a_{\ell+1} \leq a_{\ell} < a_{\ell}+\delta = c_{\ell} \leq  b_{\ell} \leq b_{\ell-1} = a_{\ell-1} = c_{\ell-1}
\]
and so the vector \mbc\ is decreasing.  This proves (i). 

If  $j \in \{1,\ldots, k-1\}$, then 
\[
\sum_{i=1}^j c_i = \sum_{i=1}^j a_i. 
\]
If  $j \in \{k,\ldots, \ell -1\}$, then
\[
\sum_{i=1}^j c_i = \sum_{i=1}^{k-1} a_i + (a_k-\delta) +  \sum_{i= k+1}^j a_i <  \sum_{i=1}^j a_i.  
\]
If  $j \in \{\ell,\ldots, n \}$, then
\[
\sum_{i=1}^j  c_i = \sum_{i=1}^{k-1} a_i + (a_k-\delta) +  \sum_{i= k+1}^{\ell-1} a_i 
+ (a_{\ell} + \delta) +  \sum_{i=\ell +1}^j a_i = \sum_{i=1}^j a_i. 
\]
Therefore, $\mbc \prec \mba$. 

Next we prove that $\mbb \prec \mbc$.  
If $j \in \{1,\ldots, k-1\}$, then 
\[
\sum_{i=1}^j b_i \leq \sum_{i=1}^j a_i = \sum_{i=1}^j c_i.
\]
We have 
\[
b_k \leq  a_k - \delta = c_k 
\]
\[
c_{\ell} = a_{\ell}+\delta \leq b_{\ell} 
\] 
\[
b_i = a_i = c_i \quad \text{ for $i \in \{k+1,\ldots, \ell -1\}$.}
\]  
Therefore, if  $j \in \{k,\ldots, \ell -1\}$, then
\[
\sum_{i=1}^j b_i =  
\sum_{i=1}^{k-1} b_i   + b_k   + \sum_{i=k+1}^{j} b_i  
\leq  \sum_{i=1}^{k-1} a_i  + (a_k-\delta)  + \sum_{i=k+1}^j c_i   = \sum_{i=1}^j c_i. 
\]
If  $j \in \{\ell,\ldots, n-1 \}$, then
\[
\sum_{i=1}^j  b_i \leq \sum_{i=1}^j  a_i 
= \sum_{\substack{i=1\\i \neq k,\ell}}^j a_i  + (a_k-\delta) + (a_{\ell} + \delta) 
 = \sum_{i=1}^j c_i. 
\]
Finally,
\[
\sum_{i=1}^n  b_i  = \sum_{i=1}^n  a_i = 
 \sum_{\substack{i=1\\i \neq k,\ell}}^n a_i  + (a_k-\delta) + (a_{\ell} + \delta) 
=  \sum_{i=1}^n c_i. 
\]
Therefore, $\mbb \prec \mbc$.  This proves (ii).  

Because  $a_{\ell} < c_{\ell} \leq c_k < a_k$  
and $a_i = c_i$ for all $i \neq k,\ell$, we have 
\[
d_H(\mba,\mbc) = 2.
\]
The triangle inequality for Hamming distance implies 
\[
2 < d_H(\mba,\mbb) \leq d_H(\mba,\mbc) + d_H(\mbb,\mbc) = 2 + d_H(\mbb,\mbc)
\]
and so $d_H(\mbb,\mbc) > 0$.   
Lemma~\ref{Muirhead:lemma:Hamming} implies $d_H(\mbb,\mbc) \geq 2$.

Because $a_k \neq b_k$ and $a_{\ell} \neq b_{\ell}$, we have  
\begin{align*}
r  = d_H(\mba,\mbb) - 2 & =  \card\{ i : b_i \neq a_i \text{ and } i \neq k, \ell \} \\ 
& = \card\{ i : b_i \neq c_i \text{ and } i \neq k, \ell \} . 
\end{align*}
In the vector \mbc, we have $c_k = b_k$ or $c_{\ell}  = b_{\ell}$, and so 
$d_H(\mbb,\mbc) = r$ or $r+1$.  Thus, 
 $d_H(\mbb,\mbc) \leq r+1 <  d_H(\mba,\mbb)$.  
This proves (iii).  

To prove~(iv), we observe that the vectors \mba\ and \mbc\ satisfy the conditions of Lemma~\ref{Muirhead:lemma:HardyLittlewoodPolya-T1} and so there is a 
Hardy-Littlewood-P\' olya $T$-transformation such that $T\mba = \mbc$. 
The Hardy-Littlewood-P\' olya $T$-transformations are doubly stochastic.
This completes the proof.  
\end{proof}

For example, in $\R^3_{\geq 0}$ the vectors 
\[
\mba = \bmat 8 \\ 3 \\ 1  \emat \qqand \mbb =  \bmat 6 \\ 4 \\ 2\emat 
\]
satisfy $\mbb \prec \mba$ and $d_H(\mbb, \mba) = 3$.  
The construction in the proof of Lemma~\ref{Muirhead:lemma:HardyLittlewoodPolya-Tc}
yields the vector 
\[
\mbc = \bmat 7 \\ 4 \\ 1 \emat
\]
which satisfies 
\[
\mbb \prec \mbc \prec \mba \qqand d_H(\mbb,\mbc) = d_H(\mba,\mbc) = 2. 
\]
The doubly stochastic matrices 
\[
T_1 = \bmat 4/5 & 1/5 & 0 \\
1/5 & 4/5 & 0 \\ 0 & 0 & 1 \emat 
\qqand
T_2 = \bmat 5/6 & 0 & 1/6 \\ 0 & 1 & 0 \\ 1/6 & 0 & 5/6 \emat 
\]
satisfy $T_1(\mba )= \mbc$ and $T_2 (\mbc) = \mbb$.  The doubly stochastic matrix
\[
T = T_2T_1 = \bmat
2/3 & 1/6 & 1/6 \\ 1/5 & 4/5 & 0\\ 2/15 & 1/30 & 5/6 
\emat
\]
satisfies $T\mba = \mbb$.

\bt                                     \label{Muirhead:theorem:HardyLittlewoodPolya-interpolation}
Let \mba\ and \mbb\ be distinct vectors in $\Rn$ such that $\mbb \prec \mba$.
For some positive integer  $r \leq d_H(\mba,\mbb)$,  
there is a sequence of vectors 
$\mbc_0, \mbc_1, \ldots, \mbc_r \in \Rn$ such that 
\beq                   \label{Muirhead:HardyLittlewoodPolya-interpolation-1}
\mbb = \mbc_r \prec \mbc_{r-1} \prec \cdots \prec \mbc_1 \prec \mbc_0 = \mba 
\eeq
and 
\beq                   \label{Muirhead:HardyLittlewoodPolya-interpolation-2}
d_H(\mbc_{i-1}, \mbc_i) = 2 \qquad \text{for all $i \in \{1,\ldots, r\}$.}
\eeq
\et

\begin{proof}
The proof is by induction on the Hamming distance $h = d_H(\mba,\mbb)$.  
We have $h \geq 2$ because $\mba \neq \mbb$.  
 
If $h= 2$, then conditions~\eqref{Muirhead:HardyLittlewoodPolya-interpolation-1} 
and~\eqref{Muirhead:HardyLittlewoodPolya-interpolation-2} are satisfied with $r=1$, 
$\mbc_0 = \mba$, and $\mbc_1 = \mbb$.  

Let $h > 2$, and assume that the Theorem holds if $d_H(\mba,\mbb) \leq h-1$.  
Let vectors \mba\ and \mbb\ satisfy $\mba \prec \mbb$ and $d_H(\mba,\mbb) = h$.  
Let $\mbc_0 = \mba$.  
By Lemma~\ref{Muirhead:lemma:HardyLittlewoodPolya-Tc}, 
there exists a vector $\mbc_1 \in \R^n$  such that 
\beq                   \label{Muirhead:HardyLittlewoodPolya-interpolation-3}
\mbb \prec \mbc_1 \prec \mbc_0 = \mba    \qqand   d_H(\mbc_0,\mbc_1) = 2 
\eeq
and  
\[
2 \leq d_H(\mbb,\mbc_1)   < d_H(\mba,\mbb).
\]
Thus, 
\[
d_H(\mbb,\mbc_1) \leq h-1. 
\]
Applying the induction hypothesis to the vectors $\mbb \prec \mbc_1$, 
we obtain an integer $r \leq h$ and a sequence of vectors 
$\mbc_1,\ldots, \mbc_r$ such that 
\beq                   \label{Muirhead:HardyLittlewoodPolya-interpolation-4}
\mbb = \mbc_r \prec \mbc_{r-1} \prec \cdots \prec \mbc_1 
\eeq
and $d_H(\mbc_{i-1}, \mbc_i) = 2$ for all $i \in \{2,\ldots, r\}$.
Conjoining the majorization chains~\eqref{Muirhead:HardyLittlewoodPolya-interpolation-3} 
and~\eqref{Muirhead:HardyLittlewoodPolya-interpolation-4} completes the induction.  
\end{proof}

Let $\mba$ and \mbb\ be vectors in \Rn\ such that $\mbb \prec \mba$.  
A sequence of vectors $\mbc_0, \mbc_1, \ldots, \mbc_r \in \Rn$ such that 
\[
\mbb = \mbc_r \prec \mbc_{r-1} \prec \cdots \prec \mbc_1 \prec \mbc_0 = \mba 
\]
and 
\[
d_H(\mbc_{i-1}, \mbc_i) = 2 \qquad \text{for all $i \in \{1,\ldots, r\}$}
\]
is called a \index{strict majorization chain}\index{majorization chain}\emph{strict majorization chain of length $r$} 
from \mba\ to \mbb.
For example, 
\[
\mbb =  \bmat 6 \\ 4 \\ 2 \emat \prec  \bmat 6 \\ 5 \\ 1  \emat \prec  \bmat 8 \\ 3 \\ 1 \emat = \mba 
\]
and 
\[
\mbb =  \bmat 6 \\ 4 \\ 2 \emat \prec  \bmat 6 \\ 5 \\ 1  \emat \prec  \bmat  7 \\ 4 \\ 1 \emat \prec  \bmat 8 \\ 3 \\ 1 \emat = \mba 
\]
are strict majorization chains of lengths 2 and 3, respectively.

\bt                                  \label{Muirhead:theorem:HardyLittlewoodPolya-T2}
Let \mba\ and \mbb\ be distinct vectors in $\Rn_{\geq 0}$.  
If $\mbb \prec \mba$, then there is a doubly stochastic matrix $P$ 
such that $\mbb = P(\mba)$.  Moreover, if $d_H(\mba,\mbb) = h$, 
then the matrix $P$ is a product 
of at most $h-1$ Hardy-Littlewood-P\' olya $T$-transformations.  
\et

\begin{proof}
The proof is by induction on the Hamming distance $d_H(\mba,\mbb)$.

We have  $d_H(\mba,\mbb) \geq 2$.  
If $d_H(\mba,\mbb) = 2$, then the vector \mbb\ satisfies the 
conditions of the vector \mbc\ 
in Lemma~\ref{Muirhead:lemma:HardyLittlewoodPolya-T1}, and so $\mbb = T\mba$ 
for some Hardy-Littlewood-P\'olya transformation $T$.

Let $h = d_H(\mba,\mbb) \geq 3$ and suppose the Theorem is true for  vectors 
\mba\ and \mbb\  with $d_H(\mba,\mbb) \leq h-1$.
Let $\mba,\mbb \in \Rn_{\geq 0}$ satisfy $\mbb \prec \mba$ and 
$d_H(\mba,\mbb) = h$.  By Lemma~\ref{Muirhead:lemma:HardyLittlewoodPolya-Tc}, 
there is a vector $\mbc \in \Rn_{\geq 0}$ such that 
\[
\mbb \prec \mbc  = T\mba   
\]
for some Hardy-Littlewood-P\' olya $T$-transformation, 
and 
\[
d_H(\mbb,\mbc) \leq d_H(\mbb,\mba) -1 \leq h-1.
\]
By the induction hypothesis, there is a doubly stochastic matrix $P'$ such that 
$\mbb = P' \mbc$ and $P'$ is the product of at most $h-2$ 
Hardy-Littlewood-P\' olya $T$-transformations.  
Then $P = P'T$ is a doubly stochastic matrix 
that is  the product of at most $h-1$ 
Hardy-Littlewood-P\' olya $T$-transformations, and 
$\mbb = P' \mbc = P' T\mba = P(\mba)$. 
This completes the proof. 
\end{proof}

\bt              \label{Muirhead:theorem:FundamentalMajorize} 
Let \mba\ and \mbb\ be vectors in $\R^n_{\geq 0}$.  
The vector \mba\ majorizes \mbb\ if and only if there exists a doubly stochastic 
matrix $P$ such that $P(\mba) = \mbb$.  
\et

 \begin{proof}
 This follows immediately from Theorems~\ref{Muirhead:theorem:DoublyStochasticMajorize} 
 and~\ref{Muirhead:theorem:HardyLittlewoodPolya-T2}. 
 \end{proof}

The following result is a fundamental theorem that relates majorization 
and convexity.

\bt                         \label{Muirhead:theorem:permutohedron-majorize}
Let \mba\ and \mbb\ be vectors in $\R^n_{\geq 0}$.  
Let $K(\mba)$ be the $S_n$-permutohedron  generated by \mba.  
The vector \mba\ majorizes \mbb\ if and only if $\mbb \in K(\mba)$. 
\et

\begin{proof}
Let $P_{\sigma}$ be the permutation matrix constructed from the permutation $\sigma \in S_n$. 
The permutohedron $K(\mba)$ is the convex hull of the set of vectors
$\{P_{\sigma}(\mba): \sigma \in S_n\}$.  
By Theorem~\ref{Muirhead:theorem:permutohedron-majorize}, 
if \mba\ majorizes \mbb, then there is a doubly stochastic matrix $P$ such that $P(\mba) = \mbb$.  
The Birkhoff-von Neumann theorem in convexity theory states that every doubly stochastic matrix 
  is a convex combination of permutation matrices.  
  It follows that there are permutations $\sigma_1,\ldots, \sigma_k \in S_n$ 
and positive numbers $\lambda_1,\ldots, \lambda_k$ such that 
\[
\sum_{i=1}^k \lambda_i = 1
\]
and 
\[
P = \sum_{i=1}^k \lambda_i P_{\sigma_i}.
\]
Therefore,
\[
\mbb = P(\mba) = \sum_{i=1}^k \lambda_i P_{\sigma_i}  (\mba)  
\] 
is a vector in $K(\mba)$.  

Conversely, if $\mbb \in K(\mba)$, then there exist permutations $\sigma_1,\ldots, \sigma_k \in S_n$ 
and positive numbers $\lambda_1,\ldots, \lambda_k$ such that 
\[
\sum_{i=1}^k \lambda_i = 1
\]
and 
\[
\mbb = \sum_{i=1}^k \lambda_i P_{\sigma_i} (\mba ) = P(\mba)
\]
where the matrix $P =  \sum_{i=1}^k \lambda_i P_{\sigma_i}$ is doubly stochastic.  
Theorem~\ref{Muirhead:theorem:DoublyStochasticMajorize} 
implies that $\mba$ majorizes \mbb.  This completes the proof. 
\end{proof}

We also have a majorization inequality for convex functions.
Let $I \subseteq \R$ be a finite or infinite interval.  
The function $f:I \rightarrow \R$ is \emph{convex} if, for all $a_1,a_2 \in I$ and $t_1, t_2 \in [0,1]$ 
with $t_1+t_2 = 1$, 
\beq                        \label{Muirhead:convexFunction-1}
f(t_1a_1+t_2a_2) \leq t_1 f(a_1)+t_2f(a_2).
\eeq
It follows by induction on $n$ 
that if $a_1,a_2,\ldots, a_n \in I$ and $t_1, t_2,\ldots, t_n \in [0,1]$ 
with $t_1+t_2 + \cdots + t_n = 1$, then 
\beq                         \label{Muirhead:convexFunction-2}
f(t_1a_1+t_2a_2 + \cdots t_na_n ) \leq t_1 f(a_1)+t_2f(a_2) + \cdots +t_n f(a_n).
\eeq

\bt                                             \label{Muirhead:theorem:convex function}  
Let $f$ be a convex function on the interval $I$. 
Let $\mba,\mbb \in \R^n_{\geq 0}$ with $\mbb \prec \mba$.  
If $\mba = \vectorsmallan$ and $\mbb = \vectorsmallbn$ 
with $a_i, b_i \in I$ for all $i \in \{1,\ldots, n\}$, then  
\[
f(b_1) + \cdots + f(b_n) \leq f(a_1) + \cdots + f(a_n). 
\]
\et

\begin{proof}
By Theorem~\ref{Muirhead:theorem:permutohedron-majorize}, 
 there is a doubly stochastic matrix $P = \bmat p_{i,j}\emat$ 
such that $\mbb = P(\mba)$.  Equivalently, 
\[
b_i = \sum_{j=1}^n p_{i,j} a_j 
\]
for all $i \in \{1,\ldots, n\}$.  
Because $P$ is doubly stochastic, we have  
\[
\rowsum_i(P) = \sum_{j=1}^n p_{i,j} =  1
\] 
for all $i \in \{1,\ldots, n\}$ and 
\[
\colsum_j(P) = \sum_{i=1}^n p_{i,j} = 1
\]
 for all $j \in \{1,\ldots, n\}$.  
The convexity of the function $f$ implies that 
\[
f(b_i) = f\left(\sum_{j=1}^n p_{i,j} a_j  \right) \leq \sum_{j=1}^n p_{i,j}  f\left( a_j  \right) 
\]
for all $i \in \{1,\ldots, n\}$ and so 
\[
\sum_{i=1}^n f(b_i) \leq \sum_{i=1}^n \sum_{j=1}^n p_{i,j}  f\left( a_j  \right) 
=  \sum_{j =1}^n\left( \sum_{i=1}^n p_{i,j} \right) f\left( a_j  \right) 
=  \sum_{j =1}^n  f\left( a_j  \right).
\]
This completes the proof.  
\end{proof}

\subsection*{Acknowledgement} 
I thank Darij Grinberg for helpful comments on this paper.

\def\cprime{$'$} \def\cprime{$'$} \def\cprime{$'$}
\providecommand{\bysame}{\leavevmode\hbox to3em{\hrulefill}\thinspace}
\providecommand{\MR}{\relax\ifhmode\unskip\space\fi MR }
\providecommand{\MRhref}[2]{\href{http://www.ams.org/mathscinet-getitem?mr=#1}{#2}
}
\providecommand{\href}[2]{#2}

\end{document}